\documentclass[12pt]{article}

\usepackage{amsmath, epsfig, cite,colordvi,xcolor}
\usepackage{amssymb}
\usepackage{amsfonts}
\usepackage{latexsym}
\usepackage{graphicx}
\usepackage{tikz}
\usepackage{float}
\usepackage{mathrsfs}[12pt]
\usepackage{bm, amscd,  amsfonts, amssymb, cite,texdraw}
\usepackage{amsmath, epsfig, cite, color, colordvi}
\usepackage{array}
\usepackage{makecell}

\usepackage[colorlinks,urlcolor=black,linkcolor=black,anchorcolor=black,citecolor=black]{hyperref}
\newtheorem{thm}{Theorem}[section]

\newtheorem{prop}[thm]{Proposition}

\newtheorem{cor}[thm]{Corollary}

\newtheorem{lem}[thm]{Lemma}

\newcommand{\qed}{{\hfill\rule{4pt}{7pt}}}
\def\pf{\noindent {\it Proof.} }
\def\Gen{{\rm Gen}}
\def\code{{\rm code}}

\numberwithin{equation}{section}

\makeatletter \@addtoreset{equation}{section} \makeatother

\tikzstyle{every node}=[circle,inner sep=1pt,fill=white!60]
\tikzstyle{tn}=[shape=circle, draw, color=black!70]

\begin{document}
\begin{center}
{\bf \large Context-free Grammars for

Permutations and Increasing Trees}\\\ \\

William Y.C. Chen$^1$ and Amy M. Fu$^2$

$^{1}$Center for Applied Mathematics\\
Tianjin University\\
Tianjin 300072, P.R. China

$^{1,2}$Center for Combinatorics, LPMC-TJKLC\\
Nankai University, Tianjin 300071, P.R. China\\ \ \\
 Email: chenyc@tju.edu.cn, fu@nankai.edu.cn

\end{center}
\noindent{\bf Abstract.}  In this paper,
we introduce the notion of a grammatical labeling to describe a recursive process of generating combinatorial objects based on a context-free grammar. For example, by labeling the ascents and descents of a Stirling permutation, we obtain a grammar  for the second-order Eulerian polynomials.  By using the grammar for $0$-$1$-$2$ increasing trees given by Dumont, we obtain a grammatical derivation of the generating function of the  Andr\'e polynomials obtained by Foata and Sch\"utzenberger, without solving a differential equation. We also find a grammar for the number $T(n,k)$  of permutations of $[n]=\{1,2,\ldots, n\}$  with $k$ exterior peaks, which was independently discovered by Ma. We demonstrate that Gessel's formula for the generating function of $T(n,k)$ can be deduced from this grammar.
Moreover, by using grammars we show that the number of the permutations of $[n]$ with $k$ exterior peaks  equals the number of increasing trees on $[n]$ with $2k+1$ vertices of even degree. A combinatorial proof of this fact is also presented.



\section{Introduction}

A context-free grammar $G$ over an alphabet $A$ is defined as a set of substitution rules replacing a letter in $A$ by a formal function over $A$.
Chen \cite{Chen}   introduced the notion of the formal derivative of a context-free grammar, and used this approach to derive combinatorial identities including
 identities on generating functions and the Lagrange inversion formula.
The formal derivative with respect to a context-free grammar satisfies the relations
just like the derivative,
\[D(u+v)=D(u)+D(v),\]
\[D(uv)=D(u)v+uD(v).\]
So the Leibniz rule is valid,
\[D^n(uv)=\sum_{k=0}^{n}{{n \choose k}D^k(u)D^{n-k}(v)}.\]
As a consequence, we see that
\[D(w^{-1})=-w^{-2}D(w),\]
since $D(ww^{-1})=0$.

The formal derivatives are also connected with the exponential generating functions. Let \[ \Gen(w,t) = \sum_{n\geq 0}{D^n(w)\frac{t^n}{n!}}\]
 for any formal function $w$. Then we have the following relations
\begin{align}
\mbox{\Gen}'(w,t) & =\mbox{\Gen}(D(w),t),\label{productd1}\\[3pt]
\mbox{\Gen}(u+v,t) & =\mbox{\Gen}(u,t)+\mbox{\Gen}(v,t),\label{sum1}\\[3pt]
\mbox{\Gen}(uv,t) & =\mbox{\Gen}(u,t)\mbox{\Gen}(v,t),\label{product1}
\end{align}
where $u,v$ and $w$ are formal functions and $\Gen'(w,t)$ means the derivative of $\Gen(w,t)$ with respect to $t$.

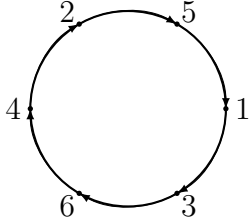
\begin{figure}[!ht]
\begin{center}
\begin{tikzpicture}
\draw[thick] (0,0) circle (1.3cm);
\coordinate [label=right:$1$] (a1) at (0:1.3cm);
\coordinate [label=60:$5$] (a2) at (60:1.3cm);
\coordinate [label=120:$2$] (a3) at (120:1.3cm);
\coordinate [label=left:$4$] (a4) at (180:1.3cm);
\coordinate [label=240:$6$] (a5) at (240:1.3cm);
\coordinate [label=300:$3$] (a6) at (300:1.3cm);
\begin{scope}[>=latex]
\draw[->] (a1) arc (0:-60:1.3cm);
\draw[->] (a6) arc (-60:-120:1.3cm);
\draw[->] (a5) arc (-120:-180:1.3cm);
\draw[->] (a4) arc (-180:-240:1.3cm);
\draw[->] (a3) arc (-240:-300:1.3cm);
\draw[->] (a2) arc (-300:-360:1.3cm);
\end{scope}
\fill (a1) circle (1pt);
\fill (a2) circle (1pt);
\fill (a3) circle (1pt);
\fill (a4) circle (1pt);
\fill (a5) circle (1pt);
\fill (a6) circle (1pt);

\end{tikzpicture}
\end{center}
\caption{a cyclic permutation of $[6]$}
\label{fig:circularper1}
\end{figure}

Dumont \cite{Dumont} introduced the following grammar
\begin{equation}\label{Dumont}
G\colon \quad x \rightarrow xy, \quad y \rightarrow xy
\end{equation}
and showed that it generates the Eulerian polynomials $A_n(x)$. For a permutation
$\pi=\pi_1\pi_2\cdots \pi_n$, the index $i\in [n-1]$ is an ascent of $\pi$ if $\pi_i<\pi_{i+1}$, a descent if $\pi_i>\pi_{i+1}$. Let $asc(\pi)$ be the number of ascents of $\pi$ and $S_n$ denote the set of permutation on $[n]=\{1, 2, \ldots, n\}$.
 The Eulerian polynomial $A_n(x)$  is defined by
\begin{equation}
 A_{n}(x)=\sum_{\pi \in S_{n}}x^{asc(\pi)+1}.
\end{equation}
To give a grammatical interpretation of $A_n(x)$,
Dumont defined bivariate polynomials $A_n(x,y)$ based on
cyclic permutations on $[n]$.
For a cyclic permutation $\sigma$, an index  $i$ $(1\leq i \leq n)$
is an ascent if $i<\sigma(i)$ and a descent if $i>\sigma(i)$. Let $asc_{c}(\sigma)$ be the number of ascents of $\sigma$, and let $des_{c}(\sigma)$ be the number of descents of $\sigma$. We assume that a cyclic permutation is oriented clockwise. For example, Figure \ref{fig:circularper1} is a cyclic permutation on $[6]$.
Let $C_n$ denote the set of cyclic permutations on $[n]$.
For $n\geq 1$, Dumont defined a the polynomial $A_n(x,y)$ as follows,
\begin{equation}\label{dumontdef1}
A_n(x,y)=\sum_{\sigma \in C_{n+1}}x^{asc_{c}(\sigma)}y^{des_{c}(\sigma)}.
\end{equation}
It should be noted that Dumont used the notation $A_{n+1}(x,y)$ instead
of $A_{n}(x,y)$ for the above polynomial. We choose the above notation
for the reason of consistency with the notation that we shall use in the
next section.

Setting $y=1$ in \eqref{dumontdef1}, we get that for $n\geq 1$
\begin{equation}\label{GenerateEulerian}
A_{n}(x,y)|_{y=1}=A_n(x).
\end{equation}
For $n\geq 1$, to obtain $A_{n+1}(x,y)$ from $A_{n}(x,y)$,  Dumont observed that the insertion of $n+1$ into a cyclic permutation of $[n]$ after $i$
 leads to a replacement of the arc $(i,\sigma(i))$ with  $(i,n+1)$   followed by
$(n+1, \sigma(i))$. If $i$ is an ascent, $(i,\sigma(i))$ corresponds to $x$ with
respect to the definition of $A_{n}(x,y)$ and the insertion of $n+1$ corresponds to
substitution of $x$ by $xy$. If $i$ is a descent, the insertion of $n+1$ corresponds to
substitution of $y$ by $xy$. Thus, $A_{n+1}(x,y)$ can be obtained from
$A_n(x,y)$ by applying the substitution rules of the grammar $G$, namely,
\begin{equation*}
A_{n+1}(x,y)=D(A_{n}(x,y)).
\end{equation*}
It follows that
\begin{equation*}
D^n(x)=A_{n}(x,y).
\end{equation*}

To demonstrate how to use a  context-free grammar to generate  combinatorial objects, we introduce the concept of a grammatical labeling. This idea is implicit in the
partition argument with respect to the grammar $f_{i}\rightarrow f_{i+1}g_1$, $ g_{i}\rightarrow g_{i+1}$ to generate partitions as given by Chen \cite{Chen}. It turns out that a grammatical labeling serves a concrete connection
between a grammar and the corresponding combinatorial structure.

This paper is organized as follows. In Section 2, we use examples to illustrate the
 notion of a grammatical labeling. We give an explanation of relation \eqref{GenerateEulerian} by labeling ascents and descents of a permutation instead of a cyclic  permutation. Similarly, by labeling ascents, descents and plateaux of a Stirling permutation, we obtain a grammatical interpretation of the second-order Eulerian polynomials.  As another example, we give a grammatical explanation of the Lah numbers by labeling the
ascents and descents of a partition into lists. We also demonstrate how to use the
formal derivative with respect to the grammar $x\rightarrow xy$, $y\rightarrow xy$ to deduce an identity on
the Eulerian polynomials.

Section 3 is devoted to the applications of the grammar $x\rightarrow xy$, $y\rightarrow x$ found by Dumont \cite{Dumont} for  the Andr\'e polynomials defined in terms of $0$-$1$-$2$ increasing trees.
As shown in Chen \cite{Chen}, a context-free grammar can be rigorously used to
derive combinatorial identities in the sense that a formal derivative
plays a role analogous to the derivative in calculus.
We shall demonstrate how to use the grammar for $0$-$1$-$2$ increasing trees given by Dumont \cite{Dumont} to give  a grammatical derivation of the generating function of the Andr\'e polynomials obtained by Foata and Sch\"utzenberger, without solving a differential equation.

In Section 4, we use the grammatical labeling to concern permutations with exterior peaks. We find that the following grammar
\begin{equation*}
G\colon \quad x \rightarrow xy, \quad  y \rightarrow x^2
\end{equation*}
can be used to generate permutations with respect to exterior peaks.
This grammar was independently discovered by Ma \cite{Ma1}.
We show that Gessel's formula for the generating function of permutations on exterior peaks can be derived by using this grammar.

In Section 5, by specializing a grammar of Dumont \cite{Dumont} for increasing trees, we find that this grammar also generates increasing trees with respect to the number of vertices with even degree. To be more specific, the degree of a vertex in
a rooted tree is meant to be the number of its children.
As a consequence, we obtain that  the number of permutations of $[n]$ with $k$ exterior peaks  equals the number of increasing trees on $[n]$ with $2k+1$ vertices of even degree.

We conclude this paper with a bijection between permutations and increasing trees which connects these two statistics.
This bijection is an extension of a correspondence between alternating permutations and  even increasing trees given by Kuznetsov, Pak and Postnikov  \cite{KPP}.

\section{Grammatical Labelings}

In order to connect a context-free grammar  to a combinatorial structure, we
associate the elements of a combinatorial structure with letters in a grammar.
 Such a labeling  scheme of a combinatorial structure is called a grammatical labeling.

For example,   consider
the following grammar given by Dumont \cite{Dumont},
\begin{equation}\label{Eulergrammar}
G\colon \quad x \rightarrow xy, \quad y \rightarrow xy.
\end{equation}
We shall use a   grammatical labeling on  permutations to show that
the Eulerian polynomial $A_n(x)$ can be expressed in terms of the
formal derivative with respect to the grammar $G$.
This labeling can be easily extended to  Stirling permutations and partitions into lists.

Denote by $A(n,m)$ the number of permutations of $[n]$ with $m-1$ ascents.
The generating function
\[
A_n(x)=\sum_{m=1}^{n} A(n,m) x^m
\]
is known as the Eulerian polynomial.

We now give a grammatical labeling on  permutations to generate the Eulerian polynomials.
Let $\pi$ be a permutation of $[n]$. An index $i$  $(1\leq i\leq n-1)$, is called an ascent if $\pi_{i}<\pi_{i+1}$, a descent if $\pi_{i}>\pi_{i+1}$.
Set $\pi_0=\pi_{n+1}=0$.
For $0 \leq i \leq n$,
if $\pi_i<\pi_{i+1}$, we label $i$ by $x$,
and
if $\pi_i>\pi_{i+1}$, we label $i$ by $y$.
With this labeling,  the weight of $\pi$ is defined as the product of the labels, that is,
\[
w(\pi)=x^{asc(\pi)+1}y^{des(\pi)+1},
\]
where $asc(\pi)$ denotes the number of ascents in $\pi$
and $des(\pi)$ denotes the number of descents in $\pi$.  As will been seen,
the polynomial  $A_{n}(x,y)$ defined in terms of the descent number and the
ascent number of a cyclic polynomial as given in \eqref{dumontdef1}
can also be expressed in terms of
the descent number and the ascent number of a permutation, that is, for $n\geq 1$,
\begin{equation*}
A_{n}(x,y)=\sum_{\pi \in S_{n}}x^{asc(\pi)+1}y^{des(\pi)+1}.
\end{equation*}

A grammatical labeling plays a role of
establishing a connection  between the action of
the formal derivative $D$
 and the insertion of the element $n+1$ into a combinatorial structure on $[n]$.
  For example, let $n=6$ and $\pi=325641$.  The grammatical labeling of
   $\pi$ is given below
\[\begin{array}{ccccccl}
&3&2&5&6&4&1\\[-5pt]
x&y&x&x&y&y&y.
\end{array}\]
If we insert $7$ after $5$, the resulting permutation and its grammatical labeling
are as follows,
\[\begin{array}{cccccccl}
 &3&2&5&7&6&4&1\\[-5pt]
x&y&x&x&y&y&y&y.
\end{array}\]
It can be seen that the insertion of $7$ after $5$ corresponds to the differentiation on the label $x$ associated with $5$.
The same argument applies to the case when the new element is inserted after an element labeled by $y$.
Hence  the action of the formal derivative $D$ on the set of weights of   permutations in $S_{n}$ gives the set of weights of permutations in $S_{n+1}$. This yields the following
grammatical expression for $A_{n}(x,y)$.

\begin{thm}\label{Grammar_P}
Let $D$ be the formal derivative with respect to  grammar \eqref{Eulergrammar}. For $n\geqslant 1$, we have
\[
D^{n}(x)=\sum_{m=1}^{n}A(n,m)x^{m}y^{n+1-m}.
\]
\end{thm}

From Theorem \ref{Grammar_P}, it  follows that $D^n(x)|_{y=1}=A_n(x)$. Here we give a grammatical proof of the following classical recurrence  for the Eulerian polynomials $A_n(x)$.

\begin{prop}For $n\geq 1$, we have
\begin{equation}\label{Eulerequation2}
A_n(x)=\sum_{k=0}^{n-1} {n \choose k} A_k(x)(x-1)^{n-1-k},
\end{equation}
where $A_0(x)=1$.
\end{prop}

\pf By the definition of grammar \eqref{Eulergrammar}, we have $D(x^{-1})=-x^{-2}D(x)=-x^{-1}y$.
Hence
\begin{equation}\label{keye2}
D(x^{-1}y)=x^{-1}D(y)+yD(x^{-1})=x^{-1}y(x-y).
\end{equation}
Since $(x-y)$ is a constant with respect to $D$, we see that
\begin{equation}\label{keye1}
D^n(x^{-1}y)=x^{-1}y (x-y)^n.
\end{equation}
By the Leibniz formula, we have for $n \geq 1$,
\begin{equation}\label{eq008}
D^n(x)=D^n(y)=D^n(x  x^{-1}y)=\sum_{k=0}^n{n \choose k} D^k(x)D^{n-k}(x^{-1}y).
\end{equation}
Substituting \eqref{keye1} into \eqref{eq008} we get
\[
(x-y)x^{-1}D^n(x)=\sum_{k=0}^{n-1} {n \choose k} x^{-1}yD^k(x)(x-y)^{n-k}.
\]
Setting $y=1$,
we arrive at  \eqref{Eulerequation2}. \qed

Next, we introduce a grammar to  generate   Stirling permutations.
Let $[n]_2$ denote the multiset $\{1^2, 2^2, \ldots, n^2\}$,
where $i^2$ stands for two occurrences of $i$.
A Stirling permutation is a permutation $\pi$ of the multiset $[n]_2$ such that  for each  $1\leq i \leq n$ the elements between
two occurrences of $i$ are larger than $i$, see Gessel and Stanley \cite{GesselStanley}. For example, $123321455664$ is a Stirling permutation
on $[6]_2$.

For a Stirling permutation $\pi=\pi_1\pi_2\cdots\pi_{2n}$, an index $i$  $(1\leq i\leq 2n-1)$, is called an ascent if $\pi_{i}<\pi_{i+1}$, a descent if $\pi_{i}>\pi_{i+1}$ and a plateaux if $\pi_{i}=\pi_{i+1}$.
We shall show that the following grammar
\begin{equation}\label{k-StirlingGr}
G\colon \quad x\rightarrow xy^2, \quad y\rightarrow xy^2
\end{equation}
can be used to generate Stirling permutations. We now give
 a grammatical labeling on   Stirling permutations.
 Let $\pi=\pi_1\pi_2\cdots\pi_{2n}$ be a Stirling permutation on $[n]_2$.
 First, we add a zero at the beginning and a zero at the end of $\pi$.
 Then we label an ascent of $0\pi_1\pi_2\cdots\pi_{2n}0$ by $x$ and label a descent or a plateau by $y$. For example, let  $\pi=244215566133$. The grammatical
  labeling of $\pi$ is given below
\[\begin{array}{ccccccccccccl}
 &2&4&4&2&1&5&5&6&6&1&3&3\\[-5pt]
x&x&y&y&y&x&y&x&y&y&x&y&y.
\end{array}\]
If we insert $77$ after the first occurrence of $4$, we get
\[\begin{array}{ccccccccccccccl}
 &2&4&7&7&4&2&1&5&5&6&6&1&3&3\\[-5pt]
x&x&x&y&y&y&y&x&y&x&y&y&x&y&y.
\end{array}\]
If we insert $77$ after the second occurrence of $1$, we get
\[\begin{array}{ccccccccccccccl}
 &2&4&4&2&1&5&5&6&6&1&7&7&3&3\\[-5pt]
x&x&y&y&y&x&y&x&y&y&x&y&y&y&y.
\end{array}\]

Notice that each Stirling permutation on $[n]_2$ can be obtained by inserting $nn$ into a Stirling permutation on $[n-1]_2$. Thus, we get a grammatical interpretation of generating function
of Stirling permutations with respect to the number of ascents.

\begin{thm}\label{StirPer}
 Let $D$ be the formal derivative with respect to grammar \eqref{k-StirlingGr}. Then we have
\begin{equation}\label{Eulerequationa}
D^n(x)=\sum_{m=1}^{n}C(n,m)x^{m}y^{2n+1-m},
\end{equation}
where $C(n,m)$ denotes the number of Stirling permutations of $[n]_2$ with $m-1$ ascents.
\end{thm}

We use the notation $C_n(x)$ as used in B\'ona \cite{Bonastir} to denote the second-order Eulerian polynomials
\[C_{n}(x)=\sum_{m=1}^{n}C(n,m)x^m.\]
From Theorem \ref{StirPer}, we see that  $D^n(x)|_{y=1}=C_{n}(x)$.

In general, we can use  the grammar
\[G\colon \quad x\rightarrow xy^r, \quad y\rightarrow xy^r\]
to generate   $r$-Stirling permutations. An $r$-Stirling permutation
is a permutation on $[n]_r=\{1^r,2^r,\ldots,n^r\}$ such that the elements between two occurrences of $i$ are not smaller than  $i$.

To conclude this section, we give the following grammar
\begin{equation}\label{LahGS}
G\colon \quad  z \rightarrow xyz, \quad x \rightarrow xy, \quad y \rightarrow xy,
\end{equation}
and we show that this grammar can be used to generate partitions of $[n]$ into lists. We call the above grammar the Lah grammar.
Recall that a partition of $[n]$ into lists is a partition of $[n]$ for which the elements of each
block are linearly ordered.
For a partition into lists, label the partition itself by $z$.
 Express a list $\sigma_1\sigma_2\cdots\sigma_m$ by $0\sigma_1\sigma_2\cdots\sigma_m0$
 and label an ascent and a descent of $0\sigma_1\sigma_2\cdots\sigma_m0$ by $x$ and $y$ respectively.
For example, let $\pi=\{325,614,7\}$.
Below is the labeling of $\pi$:
\[\begin{array}{c}
 \\[-5pt]
z
\end{array}\quad
\begin{array}{cccc}
 &3&2&5\\[-5pt]
x&y&x&y
\end{array}\quad
\begin{array}{cccc}
 &6&1&4\\[-5pt]
x&y&x&y
\end{array}\quad
\begin{array}{cc}
 &7\\[-5pt]
x&y
\end{array}
\]
Using this labeling, it can be easily seen that
 grammar \eqref{LahGS} generates   partitions into lists.

\begin{thm}\label{ThmLah} Let $C(n,k,m)$ be the number of  partitions of $[n]$ into $k$ lists with $m$ ascents.
Then, we have
\[
D^n(z)=\sum_{k=1}^{n}\sum_{m=k}^{n} C(n,k,m)x^m y^{k+n-m} z.
\]
\end{thm}

In particular, setting $y=x$, we get the grammar
\begin{equation}\label{LahS}
G\colon \quad  z \rightarrow x^2z, \quad x \rightarrow x^2,
\end{equation}
which generates the signless Lah numbers
 \[ L(n,k)={n-1 \choose k-1}\frac{n!}{k!}.\]

\begin{cor} Let $D$ be the formal derivative with respect to grammar \eqref{LahS}. Then
\begin{equation}\label{Lah}
D^n(z)=x^nz\sum_{k=1}^n L(n,k) x^{k}.
\end{equation}
\end{cor}

\section{The Andr\'e Polynomials}

In this section, we use the grammar found by Dumont \cite{Dumont}
to give  a proof of the generating function formula for the Andr\'e
polynomials without solving a differential equation.
This formula was first obtained by Foata and Sch\"utzenberger \cite{FSc}.

Recall that the Andr\'e polynomials are
defined in terms of $0$-$1$-$2$ increasing trees.
An increasing tree on $[n]$ is a rooted tree with vertex set $\{0,1,2,\ldots,n\}$ in which the labels of the vertices are increasing along any path from the root.
Note that $0$ is the root.
A $0$-$1$-$2$ increasing tree is an increasing tree in which the
 degree of any vertex is at most two. Recall that in this paper, the degree of a vertex in
a rooted tree is meant to be the number of its children.
Given a $0$-$1$-$2$ increasing tree $T$, let $l(T)$ denote the number of leaves of $T$, and $u(T)$ denote the number of  vertices of $T$ with degree $1$.
Then the Andr\'e polynomial is defined by
\[E_n(x,y)=\sum_{T}{x^{l(T)}y^{u(T)}},\]
where the sum ranges over $0$-$1$-$2$ increasing trees on $[n-1]$.

Setting $x=y=1$, $E_n(x,y)$ reduces to the $n$-th Euler number $E_{n}$, which counts both $0$-$1$-$2$ increasing trees on $[n-1]$ and alternating permutations of $[n]$, see \cite{Foata,FSc,KPP}.

Foata and Sch\"utzenberger obtained the generating function of the Andr\'e polynomials in \cite{FSc} by solving a differential equation.
Later, Foata and Han \cite{FH} found a way to compute  the generating function of $E_n(x,1)$ without solving a differential equation, or equivalently, the
generating function of $E_n(x,y)$.

Dumont \cite{Dumont} introduced the grammar
\begin{equation}\label{012}
G\colon \quad x \rightarrow xy, \quad y \rightarrow x
\end{equation}
and showed that it generates the Andr\'e polynomials $E_n(x,y)$.
This fact can be justified  intuitively in terms of the following grammatical labeling.
Given a $0$-$1$-$2$ increasing tree $T$,  a leaf of $T$ is labeled by $x$, a vertex of degree $1$ in $T$ is labeled by $y$
and a vertex of degree $2$ in $T$ is labeled by $1$.
The following figure illustrates the labeling of a $0$-$1$-$2$ increasing tree on $\{1,2,3,4,5\}$.
\begin{figure}[H]
\begin{center}
\begin{tikzpicture}
\node [tn, label=60:{$0(1)$}]{}[grow=down]
	[sibling distance=15mm,level distance=10mm]
    child {node [tn,label=180:{$2(y)$}]{}
    [sibling distance=14mm,level distance=10mm]
    	child{node [tn,label=180:{$4(x)$}]{}}
    }
    child {node [tn,label=0:{$1(1)$}]{}
    [sibling distance=9mm,level distance=10mm]
		child{node [tn,label=270:{$3(x)$}]{}}
		child{node [tn,label=270:{$5(x)$}]{}}
	};
\end{tikzpicture}
\end{center}
\caption{The labeling of a $0$-$1$-$2$ increasing tree on $\{1,2,3,4,5\}$}
\label{012tree}
\end{figure}
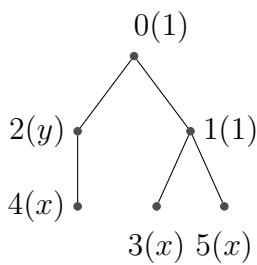
If we add $6$ as a child of $2$, the resulting tree is as follows.
 
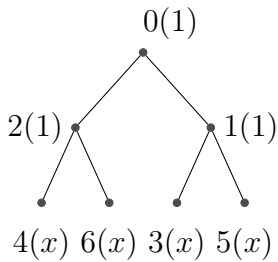
\begin{center}
\begin{tikzpicture}
\node [tn, label=60:{$0(1)$}]{}[grow=down]
	[sibling distance=18mm,level distance=10mm]
    child {node [tn,label=180:{$2(1)$}]{}
    [sibling distance=9mm,level distance=10mm]
    	child{node [tn,label=270:{$4(x)$}]{}}
        child{node[tn,label=270:{$6(x)$}]{}}
    }
    child {node [tn,label=0:{$1(1)$}]{}
    [sibling distance=9mm,level distance=10mm]
		child{node [tn,label=270:{$3(x)$}]{}}
		child{node [tn,label=270:{$5(x)$}]{}}
	};
\end{tikzpicture}
\end{center}
After the vertex $6$ is added, the label of $2$ is changed from
 $y$ to $1$, and the vertex $6$ gets a label $x$. This corresponds
 to the rule $y\rightarrow x$ of the grammar $G$.
 Similarly, adding the vertex $6$
 to a leaf of the increasing tree in Figure \ref{012tree} corresponds to the rule $x\rightarrow xy$.
 Let $D$ be the formal derivative with respect to the grammar in \eqref{012}.
The above grammatical labeling leads to the following relation
\[
D^n(x)=E_n(x,y).
\]

Now we demonstrate that one can use the grammar $G$ in \eqref{012}
to derive the generating function of $E_n(x,y)$ without solving a differential equation.

\begin{thm}[Foata and Sch\"utzenberger]\label{Eulernumgen} We have
\begin{multline}\label{newin6}
\sum_{n=0}^\infty \frac{E_n(x,y)}{n!}t^n\\[6pt]
=\frac{x\sqrt{2x-y^2}+y(2x-y^2)\sin(t\sqrt{2x-y^2})-(x-y^2)\sqrt{2x-y^2}\cos (t \sqrt{2x-y^2})}
{(x-y^2)\sin(t\sqrt{2x-y^2})+y\sqrt{2x-y^2}\cos (t \sqrt{2x-y^2})}.
\end{multline}
Setting $x=y=1$, we get
\begin{equation}\label{EulerGen}
\sum_{n=0}^\infty \frac{E_n}{n!}t^n  =\sec t +\tan t.
\end{equation}
\end{thm}


\pf
By the Leibniz rule, we have
\begin{equation}\label{eq001}
\mathrm{Gen}(x^{-1}y,t)=\mathrm{Gen}(x^{-1},t)\mathrm{Gen}(y,t).
\end{equation}
Differentiating both sides of \eqref{eq001} with respect to $t$ yields
\begin{equation}\label{eq0002}
\mathrm{Gen}'(x^{-1}y,t)
=\mathrm{Gen}'(x^{-1},t)\mathrm{Gen}(y,t)
    +\mathrm{Gen}(x^{-1},t)\mathrm{Gen}'(y,t).
\end{equation}
Since $D(x^{-1})=-x^{-1}y$, we have
\begin{equation}\label{eq0003}
\mathrm{Gen}'(x^{-1},t)=\mathrm{Gen}(D(x^{-1}),t)
=-\mathrm{Gen}(x^{-1}y,t).
\end{equation}
Using $D(y)=x$, we get
\begin{equation}\label{eq0004}
\mathrm{Gen}(x^{-1},t)\mathrm{Gen}'(y,t)=\mathrm{Gen}(x^{-1},t)\mathrm{Gen}(D(y),t)
=\mathrm{Gen}(x^{-1},t)\mathrm{Gen}(x,t)=1.
\end{equation}
Substituting \eqref{eq0003} and \eqref{eq0004} into \eqref{eq0002}, we deduce that
\[\mathrm{Gen}'(x^{-1}y,t)=1-\mathrm{Gen}(x^{-1}y,t)\mathrm{Gen}(y,t), \]
and hence
\begin{equation}\label{RecGen}
\mathrm{Gen}(y,t)=\frac{1-\mathrm{Gen}'(x^{-1}y,t)}{\mathrm{Gen}(x^{-1}y,t)}.
\end{equation}

We now compute the generating function $\Gen(x^{-1}y,t)$.
It is easily verified that
\begin{equation}\label{newin5}
D^{2m+1}(x^{-1}y)=(1-x^{-1}y^2)(y^2-2x)^m
\end{equation}
and
\begin{equation}\label{newin5-2}
D^{2m}(x^{-1}y)=x^{-1}y(y^2-2x)^m.
\end{equation}
Using \eqref{newin5} and \eqref{newin5-2}, we have
\begin{align}
\mathrm{Gen}(x^{-1}y,t)&=\sum_{n=0}^\infty \frac{D^n(x^{-1}y)}{n!}t^n \nonumber\\[6pt]
&=(1-x^{-1}y^2)\sum_{n=0}^\infty \frac{(y^2-2x)^n}{(2n+1)!} t^{2n+1}+ x^{-1}y\sum_{n=0}^\infty \frac{(y^2-2x)^n}{(2n)!}t^{2n}\nonumber\\[6pt]
&=\frac{1-x^{-1}y^2}{\sqrt{2x-y^2}}\sum_{n=0}^\infty \frac{(-1)^n(t \sqrt{2x-y^2})^{2n+1}}{(2n+1)!} + x^{-1}y\sum_{n=0}^\infty \frac{(-1)^n(t\sqrt{2x-y^2})^{2n}}{(2n)!}\nonumber\\[6pt]
&=\frac{1-x^{-1}y^2}{\sqrt{2x-y^2}}\sin(t \sqrt{2x-y^2})+x^{-1}y\cos(t\sqrt{2x-y^2}).\label{newin7}
\end{align}
Plugging \eqref{newin7} into \eqref{RecGen}, we arrive at \eqref{newin6},
and hence the proof is complete.  \qed

\section{Permutations with $k$ Exterior Peaks}

In this section, we introduce the following grammar
\begin{equation}\label{grammar}
G\colon \quad x \rightarrow xy, \quad  y \rightarrow x^2
\end{equation}
and we show that $G$ generates the number $T(n,k)$ of
permutations of $[n]$ with $k$ exterior peaks.
Let
\[
 T_n(x)=\sum_{k\geq 0}{T(n,k)x^k}.
\]
The grammar $G$ also leads to a recurrence relation of $T_n(x)$.
Moreover, we give a grammatical proof of the formula for the generating function of $T_n(x)$ due to  Gessel, see \cite{Sloane}.

Recall that for a permutation $\pi=\pi_1\pi_2\cdots \pi_n \in S_n$, the index $i$ is called an exterior peak if $1< i <n$ and $\pi_{i-1}<\pi_i>\pi_{i+1}$, or $i=1$ and $\pi_1>\pi_2$.
We shall use a grammatical labeling of permutations to show that grammar \eqref{grammar} generates the polynomial $T_n(x)$.

For a permutation $\pi$ of $[n]$, we give a labeling as follows.
First, we add an element $0$ at the end of the permutation.
If $i$ is an exterior peak, then we label $i$ and $i+1$ by $x$.
In addition, the element $0$ is labeled by $x$, and all other elements are labeled by $y$.
The weight $w$ of a permutation is defined to be the product of all the labels.
For a permutation $\pi$ with $k$ exterior peaks, its weight is given by
\[
w(\pi)=x^{2k+1}y^{n-2k}.
\]

For example, let $\pi=325641$. The labeling of $\pi$ is as follows
\[\begin{array}{ccccccl}
3&2&5&6&4&1&0\\[-5pt]
x&x&y&x&x&y&x,
\end{array}\]
and the weight of $\pi$ is $x^5 y^2$.
If we insert $7$ before $3$, then the labeling of the resulting permutation is
\[\begin{array}{cccccccl}
7&3&2&5&6&4&1&0\\[-5pt]
x&x&y&y&x&x&y&x.
\end{array}\]
We see that the label of $2$ changes from $x$ to $y$ and the label of  $7$ is $x$. So this insertion corresponds to the rule $x\rightarrow xy$.
If we insert $7$ before $0$, we get
\[\begin{array}{cccccccl}
3&2&5&6&4&1&7&0\\[-5pt]
x&x&y&x&x&y&y&x,
\end{array}\]
where the label $x$ of $0$ remains the same and the label of  $7$ is $y$.
In this case, the insertion corresponds to the rule $x\rightarrow xy$.
If we insert $7$ before $5$, we get
\[\begin{array}{cccccccl}
3&2&7&5&6&4&1&0\\[-5pt]
x&x&x&x&x&x&y&x,
\end{array}\]
where the label of $5$ changes from $y$ to $x$ and the label of  $7$ is $x$.
 This corresponds to the rule $y\rightarrow x^2$ in grammar \eqref{grammar}.
In general, the above labeling leads to the following theorem.

\begin{thm} \label{newP} Let $D$ be the formal derivative with respect to grammar \eqref{grammar}. For $n\geq 1$,
\begin{equation}\label{newp1}
D^n(x)=\sum_{k=0}^{\lfloor n/2\rfloor} T(n,k) x^{2k+1}y^{n-2k}.
\end{equation}
\end{thm}

The grammar \eqref{grammar} and relation \eqref{newp1} were announced at the International Conference on Designs, Matrices and Enumerative Combinatorics held at
National Taiwan University in 2011.
Ma \cite{Ma1} independently discovered grammar \eqref{grammar}
and gave an inductive proof of  relation \eqref{newp1}.

By Theorem \ref{newP},
we obtain the following recurrence relation.

\begin{prop} For $n\geq 1$,
\begin{equation}\label{Exeq}
T_n(x)=\sum_{j=1}^n{{n \choose j}(-1)^{j-1}(1-x)^{\lfloor j/2\rfloor}T_{n-j}(x)}.
\end{equation}
\end{prop}

\pf Note that
$$
D(x^{-1})=-x^{-1}y, \ D(-x^{-1}y)=x^{-1}(y^2-x^2), \ D(y^2-x^2)=0.
$$
Hence
\begin{equation}\label{eq0001}
D^{2m+1}(x^{-1})=-x^{-1}y(y^2-x^2)^m
\end{equation}
and
\begin{equation}\label{eq00001}
D^{2m}(x^{-1})=x^{-1}(y^2-x^2)^m.
\end{equation}
Setting $y=1$ in \eqref{eq0001} and \eqref{eq00001}, we get
\[
D^j(x^{-1})|_{y=1}=(-1)^{j}x^{-1}(1-x^2)^{{\lfloor j/2\rfloor}}.
\]
By the Leibniz rule we have
\begin{align}
D^{n}(x^{-1}x)|_{y=1}=0&=\sum_{j=0}^n {n \choose j} D^j(x^{-1})|_{y=1}D^{n-j}(x)|_{y=1}\notag\\[3pt]
&=\sum_{j=0}^n {n \choose j} (-1)^{j}x^{-1}(1-x^2)^{{\lfloor j/2\rfloor}}D^{n-j}(x)|_{y=1}.\label{RecPeak}
\end{align}
According to Theorem \ref{newP}, we see that
\[D^{n}(x)|_{y=1}=xT_n(x^2).\]
Hence  \eqref{Exeq} follows from \eqref{RecPeak}. \qed

With the aid of grammar \eqref{grammar},  we give a derivation of the following generating function of $T_n(x)$  due to Gessel, see \cite{Sloane}.

\begin{thm}[Gessel]\label{thm:Gessel} We have
\begin{equation}\label{Peak}
\sum_{n=0}^\infty \frac{T_n(x)}{n!}t^n\\
=\frac{\sqrt{1-x}}{\sqrt{1-x}\cosh(\sqrt{1-x}t)-\sinh(\sqrt{1-x}t)}.
\end{equation}
\end{thm}

To prove Theorem \ref{thm:Gessel}, we need the following generating function.
As will be seen, this generating function is related to the generating function
of $T(n,k)$.

\begin{thm}\label{lemma:Peak}
For the  the following  grammar
\begin{equation}\label{AuxiGrammar}
G\colon \quad u \rightarrow v^2, \quad v \rightarrow v,
\end{equation}
  we have
\begin{equation}\label{keyequ}
{\rm Gen}(u^{-1}v,t)
=\frac{v}{u\cosh(t)+(v^2-u)\sinh(t)}.
\end{equation}
\end{thm}

\pf Let $D$ be the formal derivative with respect to $G$.
Since $D(v)=v$, we have
\[
\Gen(v,t)=ve^t.
\]
By \eqref{product1}, we find that
\begin{equation}\label{rewe1c}
 \Gen(u^{-1}v^2,t) =\Gen(v,t)  \Gen(u^{-1}v, t)= ve^t{\rm Gen}(u^{-1}v,t) .
\end{equation}
We proceed to compute $(\Gen(u^{-1}v^2,t))'$ in two ways.
It is easily checked that
\[
D(u^{-1}v^2)=-(u^{-1}v)^2(v^2-2u).
\]
Thus, from \eqref{productd1} and \eqref{product1} we deduce that
\begin{equation}\label{rewe1a}
(\Gen(u^{-1}v^2,t))'={\rm Gen}\left(D(u^{-1}v^2),t\right)=-\Gen^2(u^{-1}v,t)\Gen(v^2-2u,t).
\end{equation}
On the other hand, since
\[
D(u^{-1}v)=u^{-1}v(1-u^{-1}v^2),
\]
from \eqref{rewe1c} we find that
\begin{align}
(\Gen(u^{-1}v^2,t))'&=(ve^t{\rm Gen}(u^{-1}v,t))'\notag\\[3pt]
&=ve^t{\rm Gen}(u^{-1}v,t)+ve^t{\rm Gen}(D(u^{-1}v),t)\notag\\[3pt]
&= ve^t\Gen(u^{-1}v,t)+ve^t\Gen(u^{-1}v,t)\Gen(1-u^{-1}v^2,t).\label{rewe1b}
\end{align}
Comparing \eqref{rewe1a} with \eqref{rewe1b}, we obtain that
\begin{align*}
  & -\Gen^2(u^{-1}v,t)\Gen(v^2-2u,t)\\[3pt]
  &  \qquad = ve^t\Gen(u^{-1}v,t)+ve^t\Gen(u^{-1}v,t)\Gen(1-u^{-1}v^2,t),
\end{align*}
or, equivalently,
\begin{equation}\label{eq:Gen}
-\Gen(u^{-1}v,t)\Gen(v^2-2u,t) = ve^t+ve^t\Gen(1-u^{-1}v^2,t).
\end{equation}
Since $D(v^2-2u)=0$, we get
\[
\Gen(v^2-2u,t)=v^2-2u.
\]
Clearly,  $\Gen(1-u^{-1}v^2,t)=1-\Gen(u^{-1}v^2,t)$.
Thus \eqref{eq:Gen} can be simplified to
\begin{equation}\label{rewe2}
-(v^2-2u)\Gen(u^{-1}v,t)=2ve^t-ve^t\Gen(u^{-1}v^2,t).
\end{equation}
Plugging \eqref{rewe1c} into \eqref{rewe2}, we arrive at
\[
{\rm Gen}(u^{-1}v,t)=\frac{2v}{v^2e^{t}-(v^2-2u)e^{-t}},
\]
which can be written in the form of \eqref{keyequ}, and so the proof is complete. \qed

We proceed to show that $\Gen(u^{-1}v,t)$ can be used to derive
the generating function of $T_n(x)$ as given in Theorem \ref{thm:Gessel}.
To this end,  we consider the following grammar
\begin{equation}\label{PeakGrammar}
G\colon\quad x\rightarrow xy, \quad y\rightarrow wx^2.
\end{equation}
For a permutation $\pi$ of $[n]$, we give a labeling which is essentially the same
as the labeling given before.
First, add an element $0$ at the end of the permutation.
If $i$ is an exterior peak, then we label $i$ by $wx$ and $i+1$ by $x$.
In addition, the element $0$ is labeled by $x$, and all other elements are labeled by $y$.
For example, let $\pi=325641$. The labeling of $\pi$ is as follows
\[
\begin{array}{ccccccc}
3&2&5&6&4&1&0\\[-5pt]
wx&x&y&wx&x&y&x.
\end{array}\]
For the grammar in \eqref{PeakGrammar}, we have
\begin{equation}\label{newp2}
D^n(x)=\sum_{k=0}^{\lfloor n/2\rfloor} T(n,k) x^{2k+1}y^{n-2k}w^k.
\end{equation}

\noindent {\it Proof of Theorem \ref{thm:Gessel}.}
For the grammar \eqref{AuxiGrammar} in Theorem \ref{lemma:Peak}, notice the relations
\begin{align*}
D(u^{-1}v)&=u^{-1}v(1-u^{-1}v^2),\\[3pt]
D(1-u^{-1}v^2)&=(u^{-1}v)^2(v^2-2u),\\[3pt]
D(v^2-2u)&=0.
\end{align*}
Comparing the above relations with the rules of the grammar in \eqref{PeakGrammar} and making the substitutions $x=u^{-1}v,\ y=1-u^{-1}v^2,\ w=v^2-2u$, we get the
rules as in grammar \eqref{PeakGrammar}, namely, $D(x)=xy$, $D(y)=wx^2$ and $D(w)=0$.
Hence   relation \eqref{newp2}
implies that
\[
D^n(u^{-1}v)=\sum_{k=0}^{\lfloor n/2\rfloor} T(n,k) (u^{-1}v)^{2k+1}(1-u^{-1}v^2)^{n-2k}(v^2-2u)^k,
\]
that is,
\begin{equation}\label{eq007}
\Gen(u^{-1}v,t)
=\sum_{n\geq 0}{\frac{t^n}{n!}
 \sum_{k=0}^{\lfloor n/2\rfloor}{T(n,k) (u^{-1}v)^{2k+1}(1-u^{-1}v^2)^{n-2k}(v^2-2u)^k}}.
\end{equation}
Comparing \eqref{keyequ} with \eqref{eq007}, we get
\[
\sum_{n\geq 0}{\frac{t^n}{n!}
\sum_{k=0}^{\lfloor n/2\rfloor}{T(n,k) (u^{-1}v)^{2k}(1-u^{-1}v^2)^{n-2k}(v^2-2u)^k}}
 =  \frac{u}{u\cosh(t)+(v^2-u)\sinh(t)}.
 \]
Since the above relation is valid for indeterminates
 $u$ and $v$, we can set $v=\sqrt{u-1}$ to deduce the following   relation
\begin{equation}\label{eq009}
\sum_{n\geq 0}{\frac{t^nu^{-n}}{n!}\sum_{k=0}^{\lfloor n/2\rfloor}{T(n,k) (1-u^2)^k}}
=\frac{u}{u\cosh(t)-\sinh(t)}.
\end{equation}
Substituting  $t$ by $ut$ in \eqref{eq009}, we get
\begin{equation}\label{eq010}
\sum_{n\geq 0}{\frac{t^n}{n!}\sum_{k=0}^{\lfloor n/2\rfloor}{T(n,k) (1-u^2)^k}}
=\frac{u}{u\cosh(ut)-\sinh(ut)}.
\end{equation}
Finally, by setting $x=1-u^2$ in \eqref{eq010}, we reach \eqref{Peak}.
This completes the proof.
\qed

\section{Peaks in permutations and  increasing trees}

In this section, we use a grammatical approach to establish the
following theorem on a
connection between permutations with a given number of exterior peaks and increasing trees with a given number of vertices of even degree.  Then we give a combinatorial
interpretation of this fact.

\begin{thm}\label{peakt}
The number of permutations $\sigma$ of $[n]$ with $k$ exterior peaks equals the number of increasing  trees $T_\sigma$ on $[n]$ with $2k+1$ vertices which have even degree.
\end{thm}

To prove the above theorem by using grammars, we first recall a grammar
 given by Dumont \cite{Dumont},
\begin{equation}\label{Euleriang}
G\colon \quad x_i \rightarrow x_0 x_{i+1}.
\end{equation}
Let $D$ be the formal derivative with respect to $G$.
Dumont \cite{Dumont} showed that
\begin{equation}\label{IncreasingTree}
D^n(x_0)=\sum_{T}{x_0^{m_0(T)}x_1^{m_1(T)}x_2^{m_2(T)}\cdots},
\end{equation}
where the sum ranges over increasing trees $T$ on $[n]$ and $m_i(T)$ denotes the number of vertices of degree $i$ in $T$.

Relation \eqref{IncreasingTree} can be justified by labeling a vertex of degree $i$ with $x_i$ in an increasing tree. Here is an example.
\begin{figure}[H]
 \begin{center}
\begin{tikzpicture}
\node [tn, label=60:{$0(x_2)$}]{}[grow=down]
	[sibling distance=19mm,level distance=10mm]
    child {node [tn,label=180:{$2(x_1)$}]{}
    [sibling distance=14mm,level distance=10mm]
    	child{node [tn,label=180:{$4(x_0)$}]{}}
    }
    child {node [tn,label=0:{$1(x_3)$}]{}
    [sibling distance=9mm,level distance=10mm]
		child{node [tn,label=270:{$3(x_0)$}]{}}
		child{node [tn,label=270:{$5(x_0)$}]{}}
        child{node [tn,label=270:{$6(x_0)$}]{}}
	};
\end{tikzpicture}
\end{center}
\caption{A labeling on an increasing tree}\label{Fig:IncreasingTree}
\end{figure}
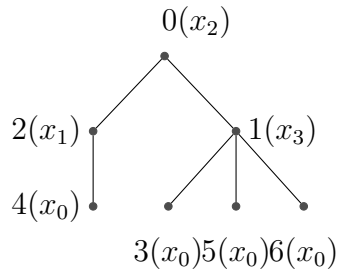

 Let $T$ be an increasing tree on $[n]$ with the above labeling. When adding the vertex $n+1$ to $T$ as the child of a vertex $v$ of degree $i$, the label of  $v$ changes from $x_i$ to $x_{i+1}$ and the label of $n+1$ is $x_0$.
This corresponds to the rule $x_i \rightarrow x_0 x_{i+1}$. Since the increasing trees on $[n+1]$ can be generated by adding $n+1$ to the increasing trees on $[n]$, the above labeling leads to \eqref{IncreasingTree}.

By setting $x_{2i}=x$ and $x_{2i+1}=y$, we see that the grammar \eqref{Euleriang} becomes the
grammar \eqref{grammar} that generates the polynomial $T_n(x)$ for permutations with
a given number of exterior peaks.
Intuitively, this leads to a
 grammatical reasoning of Theorem \ref{peakt}.
Next we give a rigorous proof of this
 observation by a grammatical labeling with respect to the parities of the
 vertices in an increasing tree.

\noindent{\it Grammatical Proof of Theorem \ref{peakt}.}
We give the following grammatical labeling of an increasing tree.
We label a vertex of even degree with $x$ and a vertex of odd degree with $y$.
For example, the labeling of the increasing tree in Figure \ref{Fig:IncreasingTree} is given below.
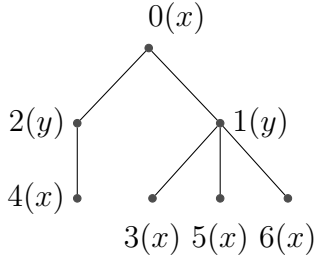
\begin{center}
\begin{tikzpicture}
\node [tn, label=60:{$0(x)$}]{}[grow=down]
	[sibling distance=19mm,level distance=10mm]
    child {node [tn,label=180:{$2(y)$}]{}
    [sibling distance=14mm,level distance=10mm]
    	child{node [tn,label=180:{$4(x)$}]{}}
    }
    child {node [tn,label=0:{$1(y)$}]{}
    [sibling distance=9mm,level distance=10mm]
		child{node [tn,label=270:{$3(x)$}]{}}
		child{node [tn,label=270:{$5(x)$}]{}}
        child{node [tn,label=270:{$6(x)$}]{}}
	};
\end{tikzpicture}
\end{center}

Let $T$ be an increasing tree on $[n]$ with the above labeling.
When adding the vertex $n+1$ to $T$ as a child of a vertex $v$ of even degree, the label of $v$ changes from $x$ to $y$ and the label of $n+1$ is $x$.
This corresponds to the rule $x \rightarrow xy$.
Similarly, adding the vertex $n+1$ to $T$ as a child of a vertex of odd degree corresponds to the rule $y \rightarrow x^2$.
Thus, we obtain that
\begin{equation}\label{SpecialDumont}
D^n(x)=\sum_{T}{x^{m_e(T)}y^{m_o(T)}},
\end{equation}
where the sum ranges over increasing trees $T$ on $[n]$ and $m_e(T)$ denotes the number of vertices of even degree in $T$, $m_o(T)$ denotes the number of vertices of odd degree in $T$.
Comparing \eqref{newp1} with \eqref{SpecialDumont}, we deduce that
\[
\sum_{k=0}^{\lfloor n/2\rfloor} T(n,k) x^{2k+1}y^{n-2k}
=\sum_{T}{x^{m_e(T)}y^{m_o(T)}},
\]
where $T$ ranges over increasing trees on $[n]$.
This completes the proof.
\qed

To conclude this paper, we give a combinatorial proof of Theorem \ref{peakt}. More precisely, we provide a bijection $\Phi$ between permutations and increasing trees such that a permutation of $[n]$ with $k$ exterior peaks corresponds to an increasing tree on $[n]$ with $2k+1$ vertices of even degree. Recall that a permutation $\sigma=\sigma_1\sigma_2\cdots \sigma_n$ of $[n]$
is called a up-down permutation if $\sigma_1<\sigma_2>\sigma_3<\cdots$.
Similarly, $\sigma$ is called a down-up permutation if $\sigma_1>\sigma_2<\sigma_3>\cdots$. When restricted to down-up permutations, $\Phi$ reduces to the bijection between down-up permutations and even increasing trees. An even increasing tree is meant to be an increasing tree such that each vertex possibly except for the root is of even degree. Kuznetsov, Pak and Postnikov \cite{KPP} gave a bijection between up-down permutations and even increasing trees.
So our bijection can be considered as an extension of the bijection given by Kuznetsov, Pak and Postnikov, since there is an obvious correspondence between
up-down permutations and down-up permutations.

Before describing our bijection, we recall that the code of a permutation is defined as follows.
For a permutation $\sigma=\sigma_1\sigma_2\cdots \sigma_n$ of $[n]$, let $\code(\sigma)=(c_1,c_2,\ldots,c_n)$ denote the code of $\sigma$. In other words, $c_i$ is the number of elements $\sigma_j$  such that $j>i$ and $\sigma_i>\sigma_j$. It is clear that $c_n=0$.

The increasing tree $\Phi(\sigma)$ can be constructed
via $n$ steps. At each step,  a vertex is added to a forest of increasing trees.
More precisely, at the $k$-th step, we obtain a forest of increasing trees with $k$ vertices, and finally obtain an increasing tree $\Phi(\sigma)$ on $[n]$.

For $k=1$, as the first step we start with an increasing tree $F_1$ with a single vertex $i_1=n-c_1$. For $k> 1$, we assume that a forest $F_{k-1}$
has been obtained at the $(k-1)$-th step. Denote by $I_{k-1}$ and $J_{k-1}$ the set of vertices and the set of roots of $F_{k-1}$. Let $\bar{I}_{k-1}$ be the complement of $I_{k-1}$, that is,
$\bar{I}_{k-1}=[n]\setminus I_{k-1}$. The goal of the $k$-th step is  to construct a  forest $F_k$
by adding an element from $\bar{I}_{k-1}$   to $F_{k-1}$.

Let $j_1,j_2,\ldots,j_l$ be the elements of $J_{k-1}$ listed in decreasing order. For notational convenience, we assume that $j_0=n+1$, $j_{l+1}=0$ and $c_0=0$.
Let
\begin{align}
U_{k}&=\{m \in \bar{I}_{k-1} \mid j_{2p+2}<m<j_{2p+1}\mbox{ for some }p\geq 0\},\label{def:U_k}\\[3pt]
V_{k}&=\{m \in \bar{I}_{k-1} \mid j_{2p+1}<m<j_{2p}\mbox{ for some }p\geq 0\}.\label{def:V_k}
\end{align}
It is clear that $U_k\cap V_k=\emptyset$ and
$ \bar{I}_{k-1}=U_k \cup V_k$.

Define $M_k$ to be $U_k$ if
$c_{k-2}\leq c_{k-1} \leq c_{k}$ or $c_{k-2}>c_{k-1} >c_{k}$;
otherwise, define $M_k$ to be $V_k$.
Let $m_1,\ldots,m_s$ be the elements of $M_k$ listed in increasing order.
We define $i_k$  to be $m_{c_{k}+1}$ if $c_{k-1} >c_{k}$, or $m_{n-k+1-c_{k}}$ if $c_{k-1} \leq c_{k}$. By the following lemma, it can be seen that
 it is feasible to choose such $i_k$, that is,
\begin{equation}\label{eq:well-definedness}
|M_{k}|\geq c_{k}+1
\end{equation}
holds if $c_{k-1}>c_{k}$ and
\begin{equation}\label{eq:well-definedness2}
|M_{k}|\geq n-k+1-c_{k}
\end{equation}
holds if $c_{k-1}\leq c_{k}$. Now, we add  $i_{k}$ to $F_{k-1}$ by  setting each $j_p\in J_{k-1}$ to be a child of $i_k$ if $j_p>i_{k}$, and let the resulting forest by $F_k$.

When $k<n$, we may iterate the above process  until we
obtain a forest $F_n$ on $[n]$. Setting each root of $F_n$
 to be a child of the vertex $0$, we obtain
 an increasing tree $T$, which is set to be  $\Phi(\sigma)$.

Here is an example for the above bijection.
Let $\sigma=5346721$.
The code of $\sigma$ is $\code(\sigma)=(4,2,2,2,2,1,0)$.
The increasing tree $\Phi(\sigma)$ is  as follows.

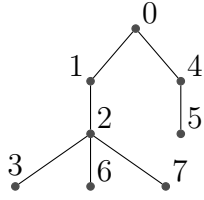
\begin{center}
\begin{tikzpicture}
\node [tn,label=60:{$0$}]{}[grow=down]
	[sibling distance=12mm,level distance=7mm]
    child{node [tn,label=120:{$1$}]{}
		[sibling distance=10mm]
    	child {node [tn,label=60:{$2$}]{}
    		child {node [tn,label=90:{$3$}]{}}
    		child {node [tn,label=60:{$6$}]{}}
			child{node [tn,label=60:{$7$}]{}}
		}
	}
	child{ node [tn,label=60:{$4$}]{}
		[sibling distance=10mm]
		child{node [tn,label=60:{$5$}]{}}
	}
;
\end{tikzpicture}
\end{center}

The values $I_k$, $J_k$, $M_{k}$, $i_k$ and the forests $F_k$ are given
in the following table.

\begin{center}
\begin{tabular}{|c|l|c|m{3cm}<{\centering}|l|}
  \hline

 $k$& \makecell[c]{$M_k$}&$i_k$ &$F_k$ &\makecell[c]{$I_k$, $J_k$}  \\ \hline
  $1$ &\makecell[c]{---} &$i_1=3$& \begin{tikzpicture}[level distance=12mm]
\node [tn,label=90:{$3$}]{};
\end{tikzpicture}&$\begin{array}{l} I_1=\{3\} \\[3pt] J_1=\{3\}\end{array}$  \\ \hline

$2$ &$M_2=\{4,5,6,7\}$&$i_2=6$& \begin{tikzpicture}[level distance=12mm]
\node [tn,label=90:{$3$}]{};
\node at(1,0) [tn,label=90:{$6$}]{};
\end{tikzpicture} &$\begin{array}{l}I_2=\{3,6\} \\[3pt] J_2=\{3,6\}\end{array}$   \\\hline

$3$ & $M_3=\{1,2,7\}$  & $i_3=7$&\begin{tikzpicture}[level distance=12mm]
\node [tn,label=90:{$3$}]{};
\node at(1,0) [tn,label=90:{$6$}]{};
\node at(2,0) [tn,label=90:{$7$}]{};
\end{tikzpicture} & $\begin{array}{l} I_3=\{3,6,7\} \\[3pt] J_3=\{3,6,7\}\end{array}$ \\\hline

$4$ & $M_4=\{1,2\}$& $i_4=2$& \begin{tikzpicture}[level distance=12mm]
\node [tn,label=90:{$2$}]{}[grow=down]
	[sibling distance=10mm,level distance=7mm]
    child {node [tn,label=270:{$3$}]{}}
    child {node [tn,label=270:{$6$}]{}}
	child{node [tn,label=270:{$7$}]{}};
\end{tikzpicture}&$\begin{array}{l}I_4=\{2,3,6,7\}\\[3pt] J_4=\{2\}\end{array}$  \\\hline

$5$ & $M_5=\{1\}$& $i_5=1$&\begin{tikzpicture}[level distance=12mm]
\node [tn,label=30:{1}]{}[grow=down]
	[sibling distance=10mm,level distance=7mm]
    child {node [tn,label=30:{$2$}]{}
    child {node [tn,label=270:{$3$}]{}}
    child {node [tn,label=270:{$6$}]{}}
	child{node [tn,label=270:{$7$}]{}}
	};
\end{tikzpicture} &$\begin{array}{l}I_5=\{1,2,3,6,7\}\\[3pt]J_5=\{1\}\end{array}$ \\\hline
$6$ &$M_6=\{4,5\}$ &$i_6=5$ &\begin{tikzpicture}[level distance=12mm]
\node [tn,label=60:{$1$}]{}[grow=down]
	[sibling distance=10mm,level distance=7mm]
    child {node [tn,label=60:{$2$}]{}
    child {node [tn,label=90:{$3$}]{}}
    child {node [tn,label=60:{$6$}]{}}
	child{node [tn,label=60:{$7$}]{}}
	};
\path(1.5,0) node [tn,label=60:{$5$}]{};
\end{tikzpicture} &$\begin{array}{l}I_6=\{1,2,3,5,6,7\} \\[3pt] J_6=\{1,5\}\end{array}$ \\\hline
$7$ &$M_7=\{4\}$ &$i_7=4$ & \begin{tikzpicture}
\node [tn,label=60:{$1$}]{}[grow=down]
	[sibling distance=10mm,level distance=7mm]
    child {node [tn,label=60:{$2$}]{}
    child {node [tn,label=90:{$3$}]{}}
    child {node [tn,label=60:{$6$}]{}}
	child{node [tn,label=60:{$7$}]{}}
	};
\path(1.5,0) node [tn,label=60:{$4$}]{}
	[sibling distance=10mm,level distance=7mm]
	child{node [tn,label=60:{$5$}]{}}
;
\end{tikzpicture}&\makecell[c]{---} \\\hline
\end{tabular}
\end{center}

In the construction of $\Phi$, at the $k$-th step ($2\leq k \leq n$)
  conditions \eqref{eq:well-definedness} and \eqref{eq:well-definedness2}
are needed to ensure the existence of the element $i_k$. The following
property implies conditions \eqref{eq:well-definedness} and \eqref{eq:well-definedness2}.

\begin{lem}\label{lem:wdef2}
For $2\leq k \leq n$, at the $k$-th step of the construction of $\Phi$, if $c_{k-1}>c_{k}$, we have that
\begin{equation}\label{well-definedness}
|M_{k}|=c_{k-1},
\end{equation}
and if $c_{k-1} \leq c_{k}$,
\begin{equation}\label{well-definedness2}
|M_{k}|=n-k+1-c_{k-1}.
\end{equation}
\end{lem}

The proof of the above lemma is parallel to the construction of $\Phi$.
First, we show that Lemma \ref{lem:wdef2} holds for $k=2$. When Lemma \ref{lem:wdef2} holds for $k$, where $k<n$,
then \eqref{eq:well-definedness} and \eqref{eq:well-definedness2} are valid for $k$, so the construction of $\Phi$
goes to the next step.

\pf
For $2\leq k\leq n$, we proceed to prove  \eqref{well-definedness} and \eqref{well-definedness2}  step by step.
It is clear that $|M_{2}|=c_{1}$ if $c_1>c_2$, and $|M_{2}|=n-1-c_{1}$ if $c_1\leq c_2$.
In other words, \eqref{well-definedness} and \eqref{well-definedness2} hold for $k=2$.
 Assume that \eqref{well-definedness} and \eqref{well-definedness2} hold for $k$.
To compute $|M_{k+1}|$, we consider the following four cases:

\noindent Case 1: $c_{k-2}>c_{k-1}>c_{k}$.
Let $j_1,j_2,\ldots,j_l$ be the elements of $J_{k-1}$ listed in decreasing order,
and let $j_0=n+1$ and $j_{l+1}=0$.
By the assumption $c_{k-2}>c_{k-1}>c_{k}$ and the definition of $M_k$,
we get
\[M_k=U_k=\{m \in \bar{I}_{k-1}\mid j_{2p+2}<m<j_{2p+1}\mbox{ for some }p\geq 0\}.\]
Since $i_k\in M_k$, there exists $q\geq 0$ such that
$j_{2q+2}<i_k<j_{2q+1}$.
So the set of roots of $F_k$ is given by
\[
J_k=\{i_k,j_{2q+2},\ldots,j_l\}.
\]
It follows that
\begin{align*}
U_{k+1}&=\{m \in \bar{I}_{k}\mid m<i_k,\,j_{2p+2}<m<j_{2p+1}\mbox{ for some }p\geq q\}\\[3pt]
&=\{m\in M_{k}\mid m<i_k\}.
\end{align*}
Since $c_{k-1}>c_{k}$, $i_k$ is the $(c_k +1)$-th smallest element in $M_k$. Hence
\[|U_{k+1}|=|\{m\in M_{k}\mid m<i_k\}|=c_k.\]
If $c_{k}>c_{k+1}$, by the assumption $c_{k-1}>c_{k}$, we have
\[
|M_{k+1}|=|U_{k+1}|=c_k.
\]
If $c_{k}\leq c_{k+1}$, we get
\[
M_{k+1}=V_{k+1}=\bar{I}_{k}\setminus U_{k+1},
\]
which implies that
\[|M_{k+1}|=n-k-c_k.\]
So we have verified that in this case \eqref{well-definedness} and \eqref{well-definedness2} are also valid for $M_{k+1}$.

\noindent Case 2: $c_{k-2}>c_{k-1}\leq c_{k}$. In this case, we have
\begin{equation}\label{mkv1}
M_k=V_k=\{m \in \bar{I}_{k-1}\mid j_{2p+1}<m<j_{2p} \mbox{ for some }p\geq 0\}.
\end{equation}
Since $i_{k}\in M_{k}$, there exists $q\geq 0$ such that $j_{2q+1}<i_k<j_{2q}$.
It follows that
\[
J_k=\{i_k,j_{2q+1},\ldots,j_l\},
\]
and hence
\[
V_{k+1}=\{m \in \bar{I}_{k}\mid i_k<m<j_0\mbox{ or }j_{2p+2}<m<j_{2p+1}\mbox{ for some }p\geq q\}.
\]
Consequently,
\begin{align}
V_{k+1}  = & \{m \in \bar{I}_{k-1}\mid j_{2p+2}<m<j_{2p+1}\mbox{ for some }p\geq 0\}\notag\\[3pt]
      & \quad +\{m \in \bar{I}_{k}\mid m>i_k,\,j_{2p+1}<m<j_{2p}\mbox{ for some }p\leq q\}. \label{V_k+1}
\end{align}
Notice that the first subset on the right hand side of \eqref{V_k+1} is
exactly $U_k$ as defined by \eqref{def:U_k}.
To compute the cardinality of the second subset on the right
hand side of \eqref{V_k+1}, we observe that $j_{2q+1}<i_k<j_{2q}$.
Hence we have
\begin{align*}
&\{m \in \bar{I}_{k}\mid m>i_k,\,j_{2p+1}<m<j_{2p}\mbox{ for some }p\leq q\}\\[3pt]
&\quad=\{m \in \bar{I}_{k}\mid m>i_k,\,j_{2p+1}<m<j_{2p}\mbox{ for some }p\geq 0\}\\[3pt]
&\quad=\{m \in M_k\mid m>i_k\}.
\end{align*}
So we obtain that
\begin{equation}\label{12vkj1}
|V_{k+1}|=|U_k|+|\{m \in M_k\mid m>i_k\}|.
\end{equation}
Using the hypothesis and \eqref{mkv1},
we find that
\begin{equation}\label{eq:V_k}
|V_k|=|M_k|=n-k+1-c_{k-1},
\end{equation}
so that
\begin{equation}\label{eq:U_k}
|U_k|=|\bar{I}_{k-1}\setminus V_k|=c_{k-1}.
\end{equation}
Since $c_{k-1}\leq c_k$, $i_k$ is the $(n-k+1-c_k)$-th smallest element in $M_k$, which implies that
\begin{equation}\label{eq:i_k}
|\{m \in M_k\mid m\leq i_k\}|=n-k+1-c_k.
\end{equation}
From \eqref{eq:V_k} and \eqref{eq:i_k} we obtain that
\begin{align}
&|\{m \in M_k\mid m>i_k\}|\notag\\[3pt]
&\quad=|M_k|-|\{m \in M_k\mid m\leq i_k\}|\notag\\[3pt]
&\quad=(n-k+1-c_{k-1})-(n-k+1-c_k)\notag\\[3pt]
&\quad=c_k-c_{k-1}.\label{eq:V_k+1-b}
\end{align}
Substituting \eqref{eq:U_k} and \eqref{eq:V_k+1-b} into \eqref{V_k+1}, we get
$|V_{k+1}|=c_k$,
and hence $|U_{k+1}|=n-k-c_k$.

If $c_{k}>c_{k+1}$, by the assumption   $c_{k-1}\leq c_{k}$, we have
$|M_{k+1}|=|V_{k+1}|=c_k.$
If $\sigma_{k}<\sigma_{k+1}$, we get
$|M_{k+1}|=|U_{k+1}|=n-k-c_k.$ This proves that in this case \eqref{well-definedness} and \eqref{well-definedness2} hold for $M_{k+1}$.

For the other two cases, $c_{k-2}\leq c_{k-1}\leq c_{k}$ and $c_{k-2}\leq c_{k-1}>c_{k}$, $|M_{k+1}|$ can be determined by the same argument. The details are omitted.
Thus we have shown that \eqref{well-definedness} and \eqref{well-definedness2} hold for $k+1$. Hence  \eqref{well-definedness} and \eqref{well-definedness2} hold for $2\leq k\leq n$. This completes the proof.
\qed

 We now have shown that $\Phi$  is well-defined.
To give a combinatorial proof of Theorem \ref{peakt}, we also need the following property.

\begin{lem}\label{lemma:Phi}
Let $\sigma=\sigma_1\sigma_2\cdots\sigma_n$ be a permutation of $[n]$ and $\code(\sigma)=(c_1,c_2,\ldots,c_n)$.
If $c_{n-1}=1$, then the root  of $\Phi(\sigma)$ is of even degree.
If $c_{n-1}=0$, then the root of $\Phi(\sigma)$ is of odd degree.
\end{lem}

\pf
Observe that for any rooted tree,
 there is an odd number of vertices of even degree. Clearly, for a permutation $\sigma$ on $[n]$, $c_{n-1}$ equals to $0$ or $1$.
It is easily seen that $c_{n-1}=1$ is equivalent to $\sigma_{n-1}>\sigma_{n}$ and
$c_{n-1}=0$ is equivalent to $\sigma_{n-1}<\sigma_{n}$. To prove the lemma, we proceed  to show that there are an odd number of non-rooted vertices of even degree in $\Phi(\sigma)$ if $\sigma_{n-1}<\sigma_{n}$,
whereas there are an even number of non-rooted vertices of even degree if $\sigma_{n-1}>\sigma_{n}$.

Recall that an index $2\leq k\leq n-1$ is called a valley of a permutation $\sigma=\sigma_1\sigma_2\cdots\sigma_n$ if $\sigma_{k-1}>\sigma_{k}<\sigma_{k+1}$.
It is clear that $i_1$ is a leaf of $\Phi(\sigma)$.
Moreover, by the construction of $\Phi$, for $2\leq k\leq n$, $i_{k}$ is a vertex of even degree if and only if $i_k\in V_k$.
Also by the construction of $\Phi$, it is easily seen that $i_k\in V_k$ if and only if $\sigma_{k-2}<\sigma_{k-1}>\sigma_{k}$ or $\sigma_{k-2}>\sigma_{k-1}<\sigma_{k}$.
Hence, for $2\leq k\leq n$, $i_{k}$ is a vertex of even degree if and only if $k-1$ is
either an exterior peak or a valley. From the above argument, we also see that $i_1$ does not correspond to any exterior peak or any valley.

We now consider the number of exterior peaks and the number of valleys in $\sigma$.
Since $\sigma_0=0$, the elements of $\sigma$ go up from $\sigma_0$, then go down
to certain position, and go up, and so on. In other words, $\sigma$ begins with an
exterior peak, then the valleys and peaks occur alternately.
If $\sigma_{n-1}<\sigma_n$, then $\sigma$ ends up with a valley.
Therefore, the number of exterior peaks equals the number of  valleys in $\sigma$.
This implies that the total number of exterior peaks and valleys is even.
Since $i_1$ is a leaf of $\Phi(\sigma)$, we see that there are an odd number of
 non-rooted vertices in $\Phi(\sigma)$ that are  of even degree. Hence the degree of $0$ must be odd.

When $\sigma_{n-1}>\sigma_n$,  $\sigma$ ends with a peak.
In this case, the number of
 exterior peaks of $\sigma$ exceeds the number of valleys of $\sigma$ by one,
 so that the total number of exterior peaks and valleys is odd.
Thus $\Phi(\sigma)$ has an even number of non-rooted vertices of even degree, since $i_1$ is a leaf.
It follows that the degree of $0$ is even, and hence the proof is complete.
\qed

We are now ready to finish the combinatorial proof of Theorem \ref{peakt}.

\noindent
{\it Combinatorial Proof of Theorem \ref{peakt}.}
We have shown that $\Phi$ is  well-defined.
To show that $\Phi$ is a bijection, we construct the inverse map $\Psi$ of $\Phi$. Let $T$ be an increasing tree on $[n]$. Start with $T$, we construct a sequence
$(c_1,c_2,\ldots,c_n)$. Let $\sigma$ be the permutation on $[n]$ such that $\code(\sigma)=(c_1,c_2,\ldots,c_n)$. Then we define  $\Psi(T)$ to be $\sigma$.

First, let $F_n$ be the forest obtained from $T$ by deleting its root $0$. Then from $F_n$, we construct a sequence of forests $F_{n-1},\ldots,F_1$. For $k=n,n-1,\ldots,2$, $F_{k-1}$ is obtained by deleting a vertex from $F_{k}$. More precisely,
for $k=n, n-1, \ldots, 2$, let $i_{k}$ be the largest root of $F_k$, and let $F_{k-1}$ be the forest obtained from $F_k$ by deleting $i_k$. For $k=1$, let $i_{1}$ be the largest root of $F_{1}$.
For $1\leq k \leq n$, let $I_k$  denote the set of vertices in $F_k$ and let $J_k$ denote
the set of roots in $F_k$.
As before, let $\bar{I}_{k}$ denote the complement of $I_{k}$ in $[n]$.
Given $I_k$ and $J_k$, assume that $U_k$ and $V_k$ are defined the same as in \eqref{def:U_k} and \eqref{def:V_k},
namely,
\begin{align*}
U_{k}&=\{m \in \bar{I}_{k-1} \mid j_{2p+2}<m<j_{2p+1}\mbox{ for some }p\geq 0\},\\[3pt]
V_{k}&=\{m \in \bar{I}_{k-1} \mid j_{2p+1}<m<j_{2p}\mbox{ for some }p\geq 0\},
\end{align*}
where $j_1,j_2,\ldots,j_l$ are the elements of $J_{k-1}$ listed in decreasing order and $j_0=n+1$, $j_{l+1}=0$.
Note that $i_k\in \bar{I}_{k-1}$ and $\bar{I}_{k-1}$ is the disjoint union of $U_k$ and $V_k$.
If $i_k\in U_k$, we set $M_k=U_k$.
If $i_k\in V_k$, we set $M_k=V_k$.

Based on $i_k$ and $M_k$, we can determine $c_k$ for $1\leq k\leq n$.
For $k=n$, it is easily seen that $|M_n|=1$. We set $c_n=0$.
For $k=n-1$, we set
\begin{equation}\label{lastone1}
c_{n-1}=\left\{
\begin{array}{ll}
1, & \textrm{if the degree of the root $0$ in $T$ is even},\\[3pt]
0, & \textrm{if the degree of the root $0$ in $T$ is odd}.\\
\end{array}
\right.
\end{equation}
Moreover, for $k=n-2,n-3,\ldots,1$, we set
\begin{equation}\label{lastone2}
c_k=\left\{
\begin{array}{ll}
|M_{k+1}|, & \textrm{if $M_{k+2}=U_{k+2}$ and $c_{k+1}>c_{k+2}$},\\[3pt]
n-k-|M_{k+1}|, & \textrm{if $M_{k+2}=U_{k+2}$ and $c_{k+1}\leq c_{k+2}$},\\[3pt]
n-k-|M_{k+1}|,  & \textrm{if $M_{k+2}=V_{k+2}$ and $c_{k+1}>c_{k+2}$},\\[3pt]
|M_{k+1}|,   & \textrm{if $M_{k+2}=V_{k+2}$ and $c_{k+1}\leq c_{k+2}$}.
\end{array}
\right.
\end{equation}
In this way, we obtain $(c_1,c_2,\ldots ,c_n)$.
Next we aim to show that for $1\leq k\leq n$,
\begin{equation}\label{wd1}
0\leq c_k \leq n-k.
\end{equation}
Since for $2\leq k \leq n$, $i_k \in M_k$ and $M_k \subseteq \bar{I}_{k-1}$,   we have
\begin{equation}\label{bonndmk}
1\leq |M_k|\leq |\bar{I}_{k-1}|.
\end{equation}
On the other hand, by the definition of $I_{k-1}$, we find that $|\bar{I}_{k-1}|=n-k+1$. It follows  that for $2\leq k \leq n$,
\begin{equation*}
1\leq |M_k|\leq n-k+1.
\end{equation*}
Clearly, for $1\leq k \leq n-1$, $c_k$ equals to $|M_{k+1}|$ or $n-k-|M_{k+1}|$.
 Thus, for $1\leq k \leq n-1$, we have
\[0\leq c_k\leq n-k.\]
Note that $c_n=0$, and so \eqref{wd1} is proved.

Let $\sigma$ be the permutation of $[n]$ with code  $(c_1,c_2,\ldots,c_n)$.
We define $\Psi(T)$ to be $\sigma$. By Lemma \ref{lem:wdef2} and Lemma \ref{lemma:Phi}, it is straightforward to verify that every step of the construction of $\Psi$ is the inverse of the corresponding step of $\Phi$. Hence $\Phi$ is a bijection.

It remains to show that $\Phi$ maps a permutation of $[n]$ with $m$ exterior peaks to an increasing tree on $[n]$ with $2m+1$ vertices of even degree. Let $\sigma$ be a permutation on $[n]$. Recall that in the proof of Lemma \ref{lemma:Phi}, we see that $\sigma$ begins with an exterior peak, then the valleys and peaks occur alternately and each peak or valley corresponds to a vertex in $[n]$ of even degree. Suppose that $\sigma$ has $m$ exterior peaks.
We shall show that $\Phi(\sigma)$ has $2m+1$ vertices of even degree.

If $\sigma_{n-1}<\sigma_n$, there are  also $m$ valleys in $\sigma$.
These $2m$ indices correspond to $2m$ vertices in $[n]$ of
even degree. As noted in the proof of Lemma \ref{lemma:Phi},
$i_1$ does not correspond to any peak or valley of $\sigma$. On the other hand, $i_1$ is a vertex of even degree since $i_{1}$ is a leaf of $\Phi(\sigma)$.
Hence, there are $2m+1$ vertices in $[n]$ of even degree in $\Phi(\sigma)$.
By Lemma \ref{lemma:Phi}, the degree of $0$ is odd.
So there are  $2m+1$ vertices of even degree in $\Phi(\sigma)$.

If $\sigma_{n-1}>\sigma_n$, there are $m-1$ valleys in $\sigma$. These $2m-1$ indices correspond to $2m-1$ vertices in $[n]$ of
 even degree. Note that $i_1$ does not correspond to any peak or valley of $\sigma$, but $i_1$ is a vertex in $[n]$ of even degree.
Hence there are $2m$ vertices in $[n]$ of even degree in $\Phi(\sigma)$.
By Lemma \ref{lemma:Phi}, the degree of $0$ is even.
So there are  $2m+1$ vertices of even degree in $\Phi(\sigma)$.
This completes the proof. \qed

\vskip 5mm \noindent{\bf Acknowledgments.} This work was supported
by the 973 Project, the PCSIRT Project of the Ministry of Education and the National
Science Foundation of China.

\end{document}